\documentclass[12pt]{amsart}
\usepackage{amsmath}
\usepackage{amsfonts,amssymb,amsthm, txfonts, pxfonts}
\def\struckint{\mathop{%
\def\mathpalette##1##2{\mathchoice{##1\displaystyle##2}%
  {##1\textstyle##2}{##1\scriptstyle##2}{##1\scriptscriptstyle##2}}%
\mathpalette
{\vbox\bgroup\baselineskip0pt\lineskiplimit-1000pt\lineskip-1000pt
\halign\bgroup\hfill$}
{##$\hfill\cr{\intop}\cr\diagup\cr\egroup\egroup}%
}\limits}

\newtheorem{theorem}{Theorem}

\newtheorem{corollary}[theorem]{Corollary}

\newtheorem{theorem-definition}[theorem]{Theorem-Definition}

\theoremstyle{remark}

\newtheorem{remark}[theorem]{Remark}

\newcommand{\ratls}{\mathbb{Q}}
\newcommand{\reals}{\mathbb{R}}
\newcommand{\sphere}{\mathbb{S}}

\begin{document}


\title[Average Distortion]{Estimates and identities for the average distortion of a linear transformation}

\author{Igor Rivin}

\address{Department of Mathematics, Temple University, Philadelphia}

\email{rivin@math.temple.edu}

\thanks{The author would like to thank Omar Hijab for enlightening discussions, and Princeton University for its hospitality during the preparation of this paper. Information on special functions has been gleaned, with gratitude, from Eric Weisstein's MathWorld.}

\date{\today}

\keywords{Linear actions, spheres, integrals}

\subjclass{37D25;49Q15}

\begin{abstract}
Let $A$ be a linear transformation $A \colon \reals^n \rightarrow \reals^n.$ We give sharp estimates on \[
\int_{\sphere^{n=1}} \log \|A u\| d u\], We also show asymptotic results (for large $n$) and evaluate a class of integrals over the sphere, including the integral of the logarithm of absolute value of one coordinate projection.
\end{abstract}

\maketitle

\section*{Introduction}
Let $A \colon \reals^n \rightarrow \reals^n$ be a non-singular linear transformation. In this note we estimate the ``average distortion'' of $A.$ More precisely, we estimate 
\begin{equation}
\label{mainthing}
\mathcal{I}(A) = \fint_{\sphere^{n-1}} \log \|A u\| d u.
\end{equation}
We remind the reader that $\fint$ denotes the mean of the integrand with respect to the measure denoted by $d u,$ which, in our case, is the standard rotationally invariant measure on the sphere. Another way of putting it is that the measure is normalized to be a probability measure. This average is of considerable interest in dynamics (see the author's paper \cite{amie} and references cited therein). A somewhat related problem of estimating the average
\begin{equation}
\label{anotherthing}
\fint_{\sphere^{n-1}}\|A u \| du
\end{equation}
is one of the principal problems addressed in the paper \cite{el3}, and we happily carry over some of the techniques and observations to the current setting.

The main estimate is summarized in the following
\begin{theorem}
\label{bigthm}
Let the singular values of the matrix $A$ be $\sigma_1, \dotsc, \sigma_n.$ Further, let 
\[
\Xi \colon\mathbb{N} \rightarrow \reals
\]
be defined as follows:
\[
\Xi(n) = 
\begin{cases}
- \log 2 - \frac12 \sum_{k=1}^{n/2} \frac1k& \text{$n$ even,}\\
-\sum_{k=1}^{(n-1)/2} \frac{1}{2k -1} & \text{$n$ odd}.
\end{cases}
\]
Then 
\[
 -  \frac12 \log n \geq \fint_{\sphere^{n-1}} \log \|A u\| d u
 - \frac12 \log \left(\sum_{i=1}^n \sigma_i^2\right)  \geq  
\Xi(n).
\]
\end{theorem}
This will be shown by way of Theorem \ref{convexity} in Section \ref{soft} followed by the explicit computation of the lower bound in Sections \ref{integration} and \ref{logintex}.  

In Section \ref{asymptotic} we will indicate asymptotic results (a ``law of large numbers'') which indicates that the upper bound in Theorem \ref{bigthm} is a better guess for reasonably well-conditioned matrices $A.$ Note, however, that the difference between the lower and upper bounds is asymptotic to $\gamma + \log 2,$ and so the gap between the two bounds is dimension independent, which indicates that $\frac12 \log \left(\sum_{i=1}^n \sigma_i^2\right)$ is the ``right'' approximation to $\mathcal{I}(A).$

\section{A sharp inequality}
\label{soft}

Our first observation is that we can assume that the matrix $A$ (as in eq.  \eqref{mainthing}) can be assumed to be diagonal, since the average in eq. \eqref{mainthing} does not change if the matrix $A$ is replaced by the diagonal matrix of its singular values $\sigma_1, \dots, \sigma_n.$ (see \cite{el3,amie} for more discussion of this), so we can replace the integral to be estimated by 
\[
I_{\Sigma} = \fint_{\sphere^{n-1}} \log \sqrt{\sum_{i=1}^n \sigma_i^2 x_i^2} d u
= \frac12 .\fint_{\sphere^{n-1}} \log \left(\sum_{i=1}^n \sigma_i^2 x_i^2\right) d u.
\]

We will estimate the quantity
\[
J_{\Sigma} = 2 I_{\Sigma} -  \log \sum_{i=1}^n \sigma_i^2.
\]
Since $J_{\Sigma}$ is scale invariant, it will be enough to estimate it (or, what is the same, $I_{\Sigma}$)
on the simplex  $H \colon \sum_{i=1}^n \sigma_i^2. $ Let us change variables, so that $s_i = \sigma_i^2.$
By the concavity of the logarithm function, the integrand of $I_{\Sigma}$ is concave (in the $(s_1, \dots, s_n)$ variables), and hence so is $I_{\Sigma}$ itself. Since, in addition, $_{\Sigma}$ is \emph{symmetric}, it follows that the maximum of $I_{\Sigma}$ on $H$ is attained when $s_1 = \dots = s_n = 1/n,$ while the minimum is attained at (any) vertex of $H,$ for example when $s_1 = 1,$ whle $s_i = 0, \quad i \neq 1.$

We have just proved the following 
\begin{theorem}
\label{convexity}
\[
 -  \frac12 \log n \geq \frac12 .\fint_{\sphere^{n-1}} \log \left(\sum_{i=1}^n \sigma_i^2 x_i^2\right) d u
 - \frac12 \log \left(\sum_{i=1}^n \sigma_i^2\right)  \geq  
\fint_{\sphere^{n-1}} \log |x_1| d u.
\]
\end{theorem}

\section{How to integrate over the sphere}
\label{integration}

\begin{theorem}
\label{logint}
Let $f: \reals^n \rightarrow \reals$ have the property that
\begin{equation}
\label{logic}
f(a \mathbf{x}) = g(a) + f(\mathbf{x}).
\end{equation}

Let $\mu$ be the standard measure on $\sphere^{n-1}.$ Then 
\begin{multline}
\label{loginteq}
2^{n/2-1} \Gamma(n/2)\int_{\sphere^{n-1}} f(x) d \mu = \\
\int_{\reals^n} f(x_1, \dotsc, x_n) \exp\left(\frac{-x_1^2 - \dotso - x_n^2}2\right) dx_1 \dots dx_n - \omega_{n-1}\int_0^\infty g(r) r^{n-1} e^{-r^2/2} d r,
\end{multline}
where 
\[
\omega_{n-1} = \dfrac{ 2 \pi^{n/2}}{\Gamma\left(\frac{n}{2}\right)},
\]
is the area of $\sphere^{n-1}.$
\end{theorem}

\begin{proof}
Let us evaluate the first integral on the right hand side. First, let us transform to polar coordinates:
\begin{multline}
I = \int_{\reals^n} f(x_1, \dotsc, x_n) \exp\left(\frac{-x_1^2 - \dotso - x_n^2}2\right) dx_1 \dots dx_n =\\
\int_0^\infty r^{n-1} e^{-r^2/2} d r \int_{\sphere^{n-1}} f(r \mathbf{x}) d \mu = \\
\int_0^\infty r^{n-1} e^{-r^2/2} dr \int_{\sphere^{n-1}} (f(x) + g(r)) d \mu = \\
\int_0^\infty r^{n-1} e^{-r^2/2} d r \int_{\sphere^{n-1}}f(x) +
\omega_{n-1} \int_0^\infty r^{n-1} g(r) e^{-r^2/2} d r
,
\end{multline}
where we have used Eq.  \eqref{logic}, and have denoted the area of the unit sphere $\sphere^{n-1}$ by $\omega_{n-1}.$
To evaluate $\omega_{n-1}$ (and, at the same time, the first integral in the right line) let $f(\mathbf{x}) = 1.$ Then $g(r) = 0,$ and we have the equation:
\[
\int_{\reals^n}\exp\left(\frac{-x_1^2 - \dotso - x_n^2}2\right) dx_1 \dots dx_n =
\omega_{n-1} \int_0^\infty r^{n-1} e^{-r^2/2} dr.
\]
The left hand side equals
\[
\left(\int_{-\infty}^{\infty} e^{-x^2/2} d x \right)^n =  (2 \pi)^{n/2}.
\]
The integral on the right hand side can be evaluated by changing variables to $u = r^2/2,$ thus getting:
\[
\int_0^\infty r^{n-1} e^{-r^2/2} d r = 2^{(n/2-1} \int_0^\infty u^{n/2 - 1} e^{-u} d u =
2^{n/2-1} \Gamma(n/2).
\]
Finally obtaining:
\[
\omega_{n-1} =\dfrac{ 2 \pi^{n/2}}{\Gamma\left(\frac{n}{2}\right)}.
\]
\end{proof}

\section{An extended example}
\label{logintex}
Let $f(x_1, \dots, x_n) = \log(|x_1|)).$ Then $f(a \mathbf(x)) = \log(a) + f(\mathbf(x)),$
so $g(r) = \log(r).$
Let us evaluate the integrals on the right hand side of \eqref{loginteq}. First,
\begin{multline}
\int_{\reals^{n}} \log(|x_1|) e^{-(x_1^2 + \dotso + x_n^2)/2} dx_1 \dots dx_n = \\
\int_{\reals^{n-1}} e^{-(x_1^2 + \dotso + x_{n-1}^2)/2} dx_1 \dots dx_{n-1} \int_{-\infty}^\infty
\log(|x|) e^{-x^2/2} d x = \\
2^{(n+1)/2} \pi^{(n-1)/2} \int_0^\infty \log(x) e^{-x^2/2} d x = \\
-2^{n/2-1} \pi^{n/2} (\gamma + \log(2)),
\end{multline}
where $\gamma$ is Euler's constant.
The second integral is also simply evaluated in terms of special functions:
\[\int_0^{\infty} \log(r) r^{n-1} e^{-r^2/2} = 2^{n/2 - 2} \Gamma(n/2) (\log(2) + \psi(n/2)),\]
where $\psi(x)$ is the logarithmic derivative of the $\Gamma$ function.
Putting everything together, we get:
\begin{multline}
\int_{\sphere^{n-1}} f(x) d\mu=\\
 \dfrac{1}{2^{n/2-1} \Gamma\left(\frac{n}2\right) }\times \\ \left(-2^{(n/2-1)} \pi^{n/2} (\gamma + \log(2)) -
\dfrac{2 \pi^{n/2}}{\Gamma\left(\frac{n}{2}\right)} 2^{n/2 - 2} \Gamma\left(\frac{n}2\right) (\log(2) + \psi(n/2))\right) = \\
- \dfrac{\pi^{n/2}}{\Gamma\left(\frac{n}2\right)} (\gamma + \log(2)) -
\dfrac{\pi^{n/2}}{\Gamma\left(\frac{n}2\right)}(\log(2) + \psi(n/2)).
\end{multline}

It is sometimes more useful to compute the \emph{mean} of a function over the sphere, and we can do that too:

\begin{multline}
\fint_{\sphere^{n-1}} f(x) d \mu = \\ \frac12 \left(-  (\gamma + \log(2)) - (\log(2) + \psi(n/2))\right) =\\
\frac12 \left(-2 \log 2 - \gamma - \psi(n/2)\right).
\end{multline}
Lest the reader is discomfitted by the appearance of the digamma function $\psi,$ we note the following simple formula for its special values at integer and half-integer points:
\[
\psi(n/2) =
\begin{cases}
-\gamma + \sum_{i=1}^{n/2-1} \frac{1}{i}, & \text{$n$ even},\\
-\gamma -2\log(2) +2\sum_{k=1}^{(n-1)/2} \frac{1}{2k-1}, & \text{otherwise}
\end{cases}
\]
Curiously, this indicates the the mean of the $\log|x_1|$ over the unit sphere is either rational or in $\ratls[\log(2)],$ since $\gamma$ always cancels.\footnote{Good thing, since it is not known whether Euler's constant is rational.}

\section{Laws of large numbers}
\label{asymptotic}
The methods of \cite{el3} go through without change to show the following results:
\begin{theorem}
\label{lawlarge}
Let $\sigma_1, \dots, \sigma_n, \dots$ be a sequence of positive numbers such that 
\[
\lim_{n \rightarrow \infty} \frac{\sum_{i=1}^n \sigma_i^4}{\left(\sum_{i=1}^n \sigma_i^2\right)^2} = 0.
\]
Let $A_n \in GL(n, \reals)$ be a matrix with singular values $\sigma_1, \dotsc, \sigma_n.$Then
\[
\lim_{n \rightarrow \infty} \fint_{\sphere^{n-1}} \| A u \| d u - \frac12\log\left(\sum_{i=1}^n \sigma_i^2\right) + \frac12 \log n = 0.
\]
\end{theorem}
\begin{corollary}
\label{condno}
The conclusion of Theorem \ref{lawlarge} holds under the assumption that there exists a constant $c$  such that $ \sigma_i/\sigma_j < c$ for all pairs $i \neq j.$
\end{corollary}

\begin{remark}
The hypothesis of Corollary \ref{condno} says that the condition numbers of the matrices $A_n$ is uniformly bounded.
\end{remark}

\bibliographystyle{plain}
\bibliography{rivin}
\end{document}